\documentclass{elsart3-1}

\usepackage{amssymb}

\usepackage{tikz}
\usepackage{tikz}
\usetikzlibrary{decorations.markings}
\usetikzlibrary{decorations.pathreplacing}
\usetikzlibrary{shapes.geometric,positioning,calc}
\usetikzlibrary{intersections}
\usetikzlibrary{positioning}

\tikzset{arrow/.style={
        decoration={markings,
            mark= at position #1 with {\arrow{stealth}},
        },
        postaction={decorate}
    }
}

\tikzset{reversearrow/.style={
        decoration={markings,
            mark= at position #1 with {\arrow{stealth reversed}},
        },
        postaction={decorate}
    }
}

\usepackage{amssymb}

\usepackage[english]{babel}
\usepackage[utf8]{inputenc}

\newtheorem{theorem}{Theorem}[section]
\newtheorem{lemma}[theorem]{Lemma}

\newtheorem{corollary}[theorem]{Corollary}
\newtheorem{e-definition}[theorem]{Definition\rm}

\newtheorem{defin}[theorem]{Definition}

\newtheorem{remark}[theorem]{\it Remark\/}

\newtheorem{condition}{Condition}

\newtheorem{theoreme}{Th\'eor\`eme}[section]

\newtheorem{proposition}[theoreme]{Proposition}

\renewcommand{\theequation}{\arabic{equation}}
\def\og{\leavevmode\raise.3ex\hbox{$\scriptscriptstyle\langle\!\langle$~}}
\def\fg{\leavevmode\raise.3ex\hbox{~$\!\scriptscriptstyle\,\rangle\!\rangle$}}

\renewcommand{\theequation}{\arabic{equation}}

\newcommand{\Fr}{\mathcal{F}}
\newcommand{\Ideal}{\mathcal{I}}
\newcommand{\Rel}{\mathcal{R}}
\newcommand{\Mon}{\mathcal{M}}
\newcommand{\SP}{\mathcal{S}}
\newcommand{\SPM}{\Lambda}
\newcommand{\fld}{k}

\newcommand{\Galg}{\fld\Fr}
\newcommand{\Qalg}{\fld\Fr / \Ideal}

\newcommand{\Dp}{\mathrm{Dp}}

\newcommand{\Ft}{\mathrm{F}}
\newcommand{\Low}{\mathrm{L}}

\newcommand{\Gr}{\mathrm{Gr}}

\newcommand{\DSpace}[1]{\langle{#1}\rangle_d}

\newcommand{\mo}[1]{\mathcal{M}({#1})}

\newcommand{\longmo}[1]{\mathcal{M}^{\geqslant 3}({#1})}

\newcommand{\mincov}[1]{\mathsf{MinCov}({#1})}

\newcommand{\nvirt}[1]{\mathsf{NVirt}({#1})}

\newcommand{\Add}{\mathrm{Add}}
\newcommand{\GreedyAlg}{\mathrm{GreedyAlg}}

\def\og{\leavevmode\raise.3ex\hbox{$\scriptscriptstyle\langle\!\langle$~}}
\def\fg{\leavevmode\raise.3ex\hbox{~$\!\scriptscriptstyle\,\rangle\!\rangle$}}

\begin{document}

\begin{frontmatter}


\selectlanguage{english}
\title{Small Cancellation Rings}




\selectlanguage{english}
\author[a]{A.\,Atkarskaya}
\ead{atkarskaya.agatha@gmail.com}
\author[b]{A.\,Kanel-Belov}
\ead{kanelster@gmail.com }
\author[c]{E.\,Plotkin }
\ead{plotkin.evgeny@gmail.com}
\author[d]{E.\,Rips}
\ead{eliyahu.rips@mail.huji.ac.il}

\address[a]{Department of Mathematics, Bar-Ilan University, Ramat Gan 5290002,
Israel and Department of Mathematics, The Hebrew University of Jerusalem, Givat Ram, 9190401 Jerusalem, Israel}
\address[b]{Department of Mathematics, Bar-Ilan University, Ramat Gan 5290002,
Israel; Department of Discrete Mathematics, Moscow Institute of Physics and Technology, Dolgoprudnyi, Institutskiy Pereulok, 141700 Moscow Oblast, Russia; College of Mathematics and Statistics, Shenzhen University, Shenzhen 518061, China}
\address[c]{Department of Mathematics, Bar-Ilan University, Ramat Gan 5290002,
Israel}
\address[d]{Institute of Mathematics, The Hebrew University of Jerusalem, Givat Ram, 9190401 Jerusalem, Israel}

\begin{abstract}

The theory of small cancellation groups is well known. In this paper we introduce the notion of Group-like Small Cancellation Ring. This is the main result of the paper. We define this ring axiomatically, by generators and defining relations. The relations must satisfy three types of axioms.
  The major one among them is called the Small Cancellation Axiom. We show that the obtained ring is non-trivial. Moreover, we show  that this ring enjoys a global filtration that agrees with relations, find a basis of the ring as a vector space  and establish the corresponding structure theorems.
It turns out that the defined ring possesses a kind of Gr\"{o}bner basis and a greedy algorithm. Finally, this ring can be used as a first step towards the iterated small cancellation theory which hopefully plays a similar role in constructing examples of rings with exotic properties as small cancellation groups do in group theory.

\vskip 0.5\baselineskip

\noindent{\bf R\'esum\'e}
\vskip 0.5\baselineskip
\noindent La th\'eorie des groupes \`a petite simplification est bien connue.
Nous introduisons la notion d'anneau \`a petite simplification de type groupe
ce que l'on peut consid\'erer comme le r\'esultat principal de cet article.
Nous donnons la d\'efinition axiomatique d'un tel anneau par g\'en\'erateurs et relations, o\`u
les relations sont d\'etermin\'ees par trois types d'axiomes dont le principal est appel\'e axiome de petite simplification.
Nous montrons que l'anneau que l'on obtient n'est pas trivial. De plus, nous montrons que cet anneau est muni
d'une filtration globale, trouvons une base de l'anneau en tant qu'espace vectoriel et \'etablissons des th\'eor\`emes de structure correspondants.
Il s'av\`ere que l'anneau que l'on a d\'efini poss\`ede une base de type Gr\"{o}bner et l'algorithme de type Dehn.
Enfin, cet anneau peut \^etre utilis\'e comme un premier pas vers la th\'eorie it\'erative de la petite simplification
qui, nous l'esp\'erons, joue le m\^eme r\^ole dans la construction d'exemples d'anneaux ayant des propri\'et\'es exotiques que jouent les groupes \`a petite simplification dans la th\'eorie des groupes.

\end{abstract}
\end{frontmatter}


\renewcommand{\labelenumi}{(\arabic{enumi})}
\newcommand{\A}{{\mathbb A}}
\newcommand{\F}{{\mathbb F}}
\renewcommand{\P}{{\mathbb P}}
\newcommand{\Q}{{\mathbb Q}}
\newcommand{\Z}{{\mathbb Z}}
\newcommand{\ko}{{\mathcal O}}
\newcommand{\kf}{{\mathcal F}}
\newcommand{\kc}{{\mathcal C}}
\newcommand{\lra}{\longrightarrow}
\newcommand{\ch}{\operatorname{char}}
\newcommand{\ord}{\operatorname{ord}}
\newcommand{\di}{\operatorname{diag}}
\newcommand{\Imm}{\operatorname{Im}}

\def\toeq{{\stackrel{\sim}{\longrightarrow}}}
\def\into{{\hookrightarrow}}


\def\alp{{\alpha}}  \def\bet{{\beta}} \def\gam{{\gamma}}
 \def\del{{\delta}}
\def\eps{{\varepsilon}}
\def\kap{{\kappa}}                   \def\Chi{\text{X}}
\def\lam{{\lambda}}
 \def\sig{{\sigma}}  \def\vphi{{\varphi}} \def\om{{\omega}}
\def\Gam{{\Gamma}}   \def\Del{{\Delta}}
\def\Sig{{\Sigma}}   \def\Om{{\Omega}}
\def\ups{{\upsilon}}


\def\F{{\mathbb{F}}}
\def\BF{{\mathbb{F}}}
\def\BN{{\mathbb{N}}}
\def\Q{{\mathbb{Q}}}
\def\Ql{{\overline{\Q }_{\ell }}}
\def\CC{{\mathbb{C}}}
\def\R{{\mathbb R}}
\def\V{{\mathbf V}}
\def\D{{\mathbf D}}
\def\BZ{{\mathbb Z}}
\def\K{{\mathbf K}}
\def\XX{\mathbf{X}^*}
\def\xx{\mathbf{X}_*}

\def\AA{\Bbb A}
\def\BA{\mathbb A}
\def\HH{\mathbb H}
\def\PP{\Bbb P}

\def\Gm{{{\mathbb G}_{\textrm{m}}}}
\def\Gmk{{{\mathbb G}_{\textrm m,k}}}
\def\GmL{{\mathbb G_{{\textrm m},L}}}
\def\Ga{{{\mathbb G}_a}}

\def\Fb{{\overline{\F }}}
\def\Kb{{\overline K}}
\def\Yb{{\overline Y}}
\def\Xb{{\overline X}}
\def\Tb{{\overline T}}
\def\Bb{{\overline B}}
\def\Gb{{\bar{G}}}
\def\Ub{{\overline U}}
\def\Vb{{\overline V}}
\def\Hb{{\bar{H}}}
\def\kb{{\bar{k}}}

\def\Th{{\hat T}}
\def\Bh{{\hat B}}
\def\Gh{{\hat G}}

\def\cF{{\mathfrak{F}}}
\def\cC{{\mathcal C}}
\def\cU{{\mathcal U}}

\def\Xt{{\widetilde X}}
\def\Gt{{\widetilde G}}

\def\gg{{\mathfrak g}}
\def\hh{{\mathfrak h}}
\def\lie{\mathfrak a}

\def\minus{^{-1}}

\def\GL{\textrm{GL}}            \def\Stab{\textrm{Stab}}
\def\Gal{\textrm{Gal}}          \def\Aut{\textrm{Aut\,}}
\def\Lie{\textrm{Lie\,}}        \def\Ext{\textrm{Ext}}
\def\PSL{\textrm{PSL}}          \def\SL{\textrm{SL}}
\def\loc{\textrm{loc}}
\def\coker{\textrm{coker\,}}    \def\Hom{\textrm{Hom}}
\def\im{\textrm{im\,}}           \def\int{\textrm{int}}
\def\inv{\textrm{inv}}           \def\can{\textrm{can}}
\def\id{\textrm{id}}              \def\Char{\textrm{char}}
\def\Cl{\textrm{Cl}}
\def\Sz{\textrm{Sz}}
\def\ad{\textrm{ad\,}}
\def\SU{\textrm{SU}}
\def\Sp{\textrm{Sp}}
\def\PSL{\textrm{PSL}}
\def\PSU{\textrm{PSU}}
\def\rk{\textrm{rk}}
\def\PGL{\textrm{PGL}}
\def\Ker{\textrm{Ker}}
\def\Ob{\textrm{Ob}}
\def\Var{\textrm{Var}}
\def\poSet{\textrm{poSet}}
\def\Al{\textrm{Al}}
\def\Int{\textrm{Int}}
\def\Mod{\textrm{Mod}}
\def\Smg{\textrm{Smg}}
\def\ISmg{\textrm{ISmg}}
\def\Ass{\textrm{Ass}}
\def\Grp{\textrm{Grp}}
\def\Com{\textrm{Com}}
\def\rank{\textrm{rank}}
\newcommand{\Or}{\operatorname{O}}

\def\tors{_\def{\textrm{tors}}}      \def\tor{^{\textrm{tor}}}
\def\red{^{\textrm{red}}}         \def\nt{^{\textrm{ssu}}}

\def\sss{^{\textrm{ss}}}          \def\uu{^{\textrm{u}}}
\def\mm{^{\textrm{m}}}
\def\tm{^\times}                  \def\mult{^{\textrm{mult}}}

\def\uss{^{\textrm{ssu}}}         \def\ssu{^{\textrm{ssu}}}
\def\comp{_{\textrm{c}}}
\def\ab{_{\textrm{ab}}}

\def\et{_{\textrm{\'et}}}
\def\nr{_{\textrm{nr}}}

\def\nil{_{\textrm{nil}}}
\def\sol{_{\textrm{sol}}}
\def\End{\textrm{End\,}}

\def\til{\;\widetilde{}\;}


\selectlanguage{english}


\section{Introduction}\label{intro}
The Small Cancellation Theory for groups is well known (see  \cite{LS}). The similar theory exists also for semigroups and monoids (see \cite{Hig}, \cite{GS}, \cite{S}). However, the construction of such a theory for systems with two operations faces  severe difficulties.

In the present paper we develop  a small cancellation  theory for associative algebras with a basis of invertible elements.  In fact, in course of studying the question:
\medskip

\centerline{\it ``what is a small cancellation  associative ring?"}

\medskip
\noindent
 we axiomatically define  a ring, which can reasonably be called a ring with  small cancellation properties and conditions. We also determine the structure and properties of this ring.

\subsection{Motivation, objectives, results}\label{motivations}

  The motivation for developing a ring-theoretical analog of small cancellation comes from the fact that small cancellation for groups and, especially, its more far-reaching versions, provide a very powerful technique for constructing groups with unusual, and even exotic, properties, like for example,  infinite Burnside groups \cite{NA1}--\cite{NA3}, \cite{A}, \cite{Ol3}, \cite{Ivanov}, \cite{Lys},  Tarski monster \cite{Ol2}, finitely generated infinite divisible groups \cite{Gu},  and many others, see e.g., \cite{Ol1}.

On the other hand, there is a conceptual desire to understand what negative curvature could mean for ring theory.

For any group with fixed system of generators, its Cayley graph can be considered as a metric space.  This leads to Gromov's program ``Groups as geometric objects" \cite{Gr1}, see also \cite{Gr2}. In particular, a finitely generated group is word-hyperbolic when its Cayley graph is $\delta$-hyperbolic for $\delta > 0$ (see \cite{Bo}, \cite{DK} for modern exposition and references).

So far, we do not know a way to associate a geometric object to a ring. Thus, having in mind the
negative curvature as a heuristic and indirect hint for our considerations, we, nevertheless,  follow a more accessible  combinatorial line of studying rings. Therefore, small cancellation groups appear naturally at the stage.   

Finitely generated small cancellation groups turned out to be word hyperbolic (when every relation needs at least 7 pieces). So, if we could generalize small cancellation to the ring theoretic situation, it would provide examples to the yet undefined concept of a ring with a negative curvature. Another source  of potential examples are group algebras of hyperbolic groups.

Following this reasoning,  we introduce in the paper the three types of axioms for rings called Compatibility Axiom, Small Cancellation Axiom and Isolation Axiom. We study  rings $\mathcal A$ with the basis of invertible elements that satisfy these axioms with respect to a fixed natural constant $\tau \geqslant 10$. We show that:
\medskip

\begin{itemize}

\item [$\bullet$] Such rings $\mathcal A$ are non-trivial;

\item [$\bullet$] Such rings $\mathcal A$ enjoy a global filtration that agrees with the relations;

\item [$\bullet$] 
An explicit basis of $\mathcal A$ as a linear space is constructed and the corresponding structure theorems are proved;

\item [$\bullet$] These rings  possess algorithmic properties similar to the ones valid for
groups with small cancellation. In particular they have solvable equality problem and enjoy a greedy algorithm;

\item [$\bullet$] These rings also possess a Gr\"{o}bner basis with respect to some sophisticated  linear order on monomials.

\end{itemize}
\medskip
The list of facts above can be viewed as a major result of the paper.  In what follows we describe  and illuminate all these items. The detailed exposition of these results is contained in the  paper  \cite{AKPR1}. Note that the axiomatic theory presented in this paper is modeled after a particular case we have treated in \cite{AKPR}.


\subsection{Small cancellation groups, background}\label{smcgroups}
Consider a group presentation $G = \langle \mathcal{X} \mid \Rel\rangle$ where we assume that the set of relations $\Rel$ is closed under cyclic permutations and inverses and all elements of $\Rel$ are cyclically reduced. The interaction between the defining relations is described in terms of  small pieces. A word $s$ is called a {\it small piece}  with respect to $\Rel$ (in generalized group sense, see \cite{Rips}, \cite{LS}) if there are relations of the form $ s r_1$ and  $s r_2$ in $\Rel$ such that  $r_1 r_2^{-1} \neq 1$ and  $r_1 r_2^{-1}$ is not conjugate to a relator from $\Rel$ in the corresponding free group, even after possible cancellations.

\noindent
{\bf Remark. } The geometric way to think about small pieces is seeing them as words that may appear on the common boundary between two cells in the van Kampen diagram \cite{Ol1}, \cite{LS}. In particular, if $r_1r_2^{-1} \in \Rel$, then we can substitute these cells by a simple cell, so we are entitled to assume from the beginning that $r_1r_2^{-1} \notin \Rel$.

The {\it small cancellation condition} says that any relation 
in $\Rel$ cannot be written as a product of too few small pieces.  For most purposes seven small pieces suffice since  the discrete Euler characteristic per cell becomes negative \cite{LS}, \cite{Ly}.

To ensure this,  we can assume that the length of any small piece is less than one sixth of the length of the relation in which it appears. The Main Theorem of Small Cancellation Theory can be stated as follows.

Let $w_1$, $w_2$ be two words that do not contain occurrences of more than a half of a relation from $\Rel$. They represent the same element of $G$ if and only if they can be connected by a one-layer diagram (\cite{LS}, especially see Greendlinger's Lemma). The transition from $w_1$ to $w_2$ can be divided into a sequence of  elementary steps called {\it turns} \cite{NA1}-\cite{NA3}. Each turn reverses just  one cell.

\subsection{Small cancellation axioms for the ring case}\label{smcrings}

First of all, given a field $k$ and the  free group $\Fr$, denote by $\Galg$ the corresponding  group algebra. Elements of $\Fr$ and $\Galg$ are called {\it monomials or words} and {\it polynomials}, respectively.
 Let a set of polynomials $\Rel$ from $\Galg$ be fixed. Define $\Ideal$ to be the ideal generated by the elements of $\Rel$.

 Let the free group $\Fr$ be freely generated by an alphabet $S$. Assume
$$
\Rel = \left\lbrace p_i =  \sum\limits_{j = 1}^{n(i)} \alpha_{ij} m_{ij} \mid \alpha_{ij}\in k, m_{ij} \in \Fr, i \in I \right\rbrace
$$
is a (finite or infinite) set of polynomials that generates the ideal $\Ideal$ (as an ideal). We denote this way of generating by $\langle \rangle_i$. So,
$$
\Ideal = \left\langle \Rel \right\rangle_i = \left\langle p_i = \sum\limits_{j = 1}^{n(i)} \alpha_{ij} m_{ij} \mid \alpha_{ij}\in k, m_{ij} \in \Fr, i \in I \right\rangle_i.
$$
We assume that the monomials $m_{ij}$ are reduced, the polynomials $p_i$ are additively reduced, 
$I$ is some index set. In particular, we assume that all coefficients $\alpha_{ij}$ are non-zero. Denote the set of all monomials $m_{ij}$ of $\Rel$ by $\Mon$.
Throughout the paper we reserve small Greek letters for non-zero elements of the field $k$.

\begin{condition}[Compatibility Axiom]
\label{comp_ax}
The axiom consists of the following two conditions.
\begin{enumerate}
\item
If $p = \sum\limits_{j = 1}^{n} \alpha_{j} m_{j} \in \Rel$, then $\beta p = \sum\limits_{j = 1}^{n} \beta\alpha_{j} m_{j} \in \Rel$ for every $\beta \in k, \beta \neq 0$.

\item
Let $x \in S \cup S^{-1}$, $p = \sum\limits_{j = 1}^{n} \alpha_{j} m_{j} \in \Rel$. Suppose there exists $j_0 \in \lbrace 1, \ldots, n \rbrace$ such that $x^{-1}$ is the initial symbol of $m_{j_0}$. Then
$$
xp = \sum\limits_{j = 1}^{n} \alpha_{j} xm_{j} \in \Rel
$$
{\rm (}after the cancellations in the monomials $xm_{j}${\rm)}.
\end{enumerate}
We require the same condition from the right side as well.
\end{condition}
From the second condition of Compatibility Axiom it immediately follows that the set $\Mon$ is closed under taking subwords. In particular, $1$ always belongs to $\Mon$.

Now we state a definition of \emph{a small piece}. It plays a central role in the further argument.

\begin{defin}
\label{sp}
Let $c \in \Mon$. Assume there exist two polynomials
$$
p = \sum\limits_{j = 1}^{n_1} \alpha_j a_j + \alpha a \in \Rel,\\ \qquad
q = \sum\limits_{j = 1}^{n_2} \beta_j b_j + \beta b \in \Rel,
$$
such that $c$ is a subword of $a$ and a subword of $b$. Namely,
$$
a = \widehat{a}_1c\widehat{a}_2, \qquad \ b = \widehat{b}_1c\widehat{b}_2,
$$
where $\widehat{a}_1$, $\widehat{a}_2$, $\widehat{b}_1$, $\widehat{b}_2$ are allowed to be empty. Assume that
$$
\widehat{b}_1 \widehat{a}_1^{-1}p = \widehat{b}_1 \widehat{a}_1^{-1}(\sum\limits_{j = 1}^{n_1} \alpha_j a_j + \alpha \widehat{a}_1c\widehat{a}_2) = \sum\limits_{j = 1}^{n_1} \alpha_j \widehat{b}_1\widehat{a}_1^{-1}a_j + \alpha \widehat{b}_1c\widehat{a}_2 \notin \Rel
$$
{\rm (}even after the cancellations  {\rm )}, or
$$
p\widehat{a}_2^{-1}\widehat{b}_2 = (\sum\limits_{j = 1}^{n_1} \alpha_j a_j  + \alpha \widehat{a}_1c\widehat{a}_2)\widehat{a}_2^{-1}\widehat{b}_2 = \sum\limits_{j = 1}^{n_1} \alpha_j a_j\widehat{a}_2^{-1}\widehat{b}_2 + \alpha \widehat{a}_1c\widehat{b}_2 \notin \Rel
$$
(even after the cancellations). Then the monomial $c$ is called a small piece.
\end{defin}

We denote the set of all small pieces by $\SP$. Clearly, $\SP \subseteq \Mon$. From the definition it follows that the set $\SP$ is closed under taking subwords. In particular, if the set $\SP$ is non-empty, the monomial $1$ is always a small piece. If the set $\SP$ is turned out to be empty, then \emph{we still assign $1$ to be a small piece}.

Let $u \in \Mon$. Then either $u = p_1\cdots p_k$, where $p_1, \ldots, p_k$ are small pieces, or $u$ cannot be represented as a product of small pieces. We introduce a measure on monomials of $\Mon$ (aka \emph{$\SPM$-measure}). {\it We say that}
$\SPM(u) = k$ \textit{ if } $u$ \textit{ can be represented as a product of small pieces} \textit{and minimal possible number of small pieces} \textit{in such representation is equal to } $k$. {\it We say that} $\SPM(u) = \infty$ \textit{ if } $u$ \textit{ can not be represented as a product of small pieces}.

We fix a constant $\tau \in \mathbb{N}$.

\begin{condition}[Small Cancellation Axiom]
\label{sc_ax}
Assume $p_1, \ldots, p_n \in \Rel$ and a linear combination $\sum\limits_{s = 1}^{n} \gamma_s p_s$ is non-zero after additive cancellations. Then there exists a monomial $a$ in $\sum\limits_{s = 1}^{n} \gamma_s p_s$ with a non-zero coefficient after additive cancellations such that either $a$ can not be represented as a product of small pieces or every representation of $a$ as a product of small pieces contains at least $\tau + 1$ small pieces. That is, $\SPM(a) \geqslant \tau + 1$, including $\SPM(a) = \infty$.
\end{condition}

\begin{defin}
\label{incident_momom}
Let $p = \sum_{j = 1}^n \alpha_j a_j \in \Rel$. Then we call the monomials $a_{j_1}, a_{j_2}$, $1 \leqslant j_1, j_2 \leqslant n$, {incident monomials} {\rm (}including the case $a_{j_1} = a_{j_2}${\rm )}. Recall that $\alpha_j \neq 0$, $j = 1, \ldots, n$.
\end{defin}

Now we introduce the last condition, we call it Isolation Axiom. Unlike two previous axioms this is entirely a ring-theoretic condition.  Here we use the  notions of  maximal occurrence of a monomial of $\Mon$ and of overlap (see Subsection~\ref{filtration_1}).
The complexity of formulation of Isolation Axiom may perflex the reader. This axiom works in the transition from monomials to tensor products and, thus, to structure theory of rings with small cancellation. It imposes essential constraints on  rings under consideration. That is why we have chosen its weakest form to make the corresponding class of rings wider. This resulted in a somewhat cumbersome  definition.
\begin{condition}[Isolation Axiom, left-sided]
Let $m_1, m_2, \ldots, m_{k}$ be a sequence of monomials of $\Mon$ such that $m_1 \neq m_k$ and $m_i, m_{i + 1}$ are incident monomials for all $i = 1, \ldots, k - 1$, and $\SPM(m_i) \geqslant \tau - 2$ for all $i = 1, \ldots, k$. Let us take a monomial $a \in \Mon$ with the following properties.
\begin{itemize}
\item [\rm {1.}]
$\SPM(a) \geqslant \tau - 2$;
\item [\rm {2.}]
$am_1, am_k \notin \Mon$, $am_1$ has no cancellations, $am_k$ has no cancellations;
\item [\rm {3.}]
$m_1$ is a maximal occurrence in $am_1$, $m_k$ is a maximal occurrence in $am_k$.
\item [\rm {4.}]
Let $ap_1(a)$ be a maximal occurrence in $am_1$ that contains $a$, let $ap_k(a)$ be a maximal occurrence in $am_k$ that contains $a$ {\rm (}that is, $p_1(a)$ is the overlap of $ap_1(a)$ and $m_1$, $p_1(a)$ may be empty, and $p_k(a)$ is the overlap of $ap_k(a)$ and $m_k$, $p_k(a)$ may be empty{\rm)}. Assume that there exist monomials $l$, $l^{\prime} \in \Mon$ such that
\begin{itemize}
\item
$l$, $l^{\prime}$ are small pieces;
\item
$la, l^{\prime}a \in \Mon$, $la$ has no cancellations, $l^{\prime}a$ has no cancellations;
\item
there exists a sequence of monomials $b_1, \ldots, b_n$ from $\Mon$ such that $b_1 = lap_1(a)$, $b_n = l^{\prime}ap_k(a)$, $b_i, b_{i + 1}$ are incident monomials for all $i = 1, \ldots, n - 1$, and $\SPM(b_i) \geqslant \tau - 2$ for all $i = 1, \ldots, n$.
\end{itemize}
\end{itemize}
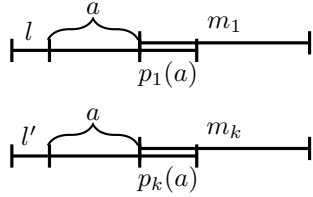
\begin{center}
\begin{tikzpicture}
\draw[|-|, black, very thick] (2.2,0.1)--(4.5,0.1) node[midway, above] {$m_1$};
\draw[|-|, black, very thick] (2.2, 0)--(3, 0) node[midway, below] {$p_1(a)$};
\draw[|-, black, very thick] (1,0)--(2.2,0);
\draw[|-, black, very thick] (0.5, 0)--(1,0) node[midway, above] {$l$};
\draw [thick, decorate, decoration={brace, amplitude=8pt, raise=4pt}] (1, 0) to node[midway, above, yshift=10pt] {$a$} (2.2, 0);
\end{tikzpicture}

\begin{tikzpicture}
\draw[|-|, black, very thick] (2.2,0.1)--(4.5,0.1) node[midway, above] {$m_k$};
\draw[|-|, black, very thick] (2.2, 0)--(3, 0) node[midway, below] {$p_k(a)$};
\draw[|-, black, very thick] (1,0)--(2.2,0);
\draw[|-, black, very thick] (0.5, 0)--(1,0) node[midway, above] {$l^{\prime}$};
\draw [thick, decorate, decoration={brace, amplitude=8pt, raise=4pt}] (1, 0) to node[midway, above, yshift=10pt] {$a$} (2.2, 0);
\end{tikzpicture}
\end{center}

Notice that since $a$ is not a small piece, then 
we get that
$lap_1(a), l^{\prime}ap_k(a) \in \Mon$, and $lap_1(a)$ is a maximal occurrence in $lap_1(a)m_1$, $l^{\prime}ap_k(a)$ is a maximal occurrence in $l^{\prime}ap_k(a)m_k$.

Then we require that ${p_1(a)}^{-1}\cdot m_1 \neq {p_k(a)}^{-1}\cdot m_k$ for every such $a \in \Mon$.
\end{condition}
The right-sided Isolation Axiom is formulated symmetrically.
\begin{remark}
We shall informally explain the essence of Isolation axioms. Given a monomial $U$, consider the set of its non-degenerate derived monomials (see Subsection 1.5). Every derived monomial can be imagined as a result of a sequence of replacements of virtual members of a chart by incident monomials. If two essentially different sequences of replacements result in one and the same derived monomial, the exotic dependencies appear in the ideal $\Ideal$. Isolation axiom guarantees that essentially different sequences of replacements result in different monomials. Hence, exotic dependencies are not present in $\Ideal$.
\end{remark}

\begin{defin}\label{smcRing}
We say that  ${\mathcal A} = \Galg/\Ideal (\Rel)$ is $C{(\tau)}$-{\it small cancellation ring} if it satisfies Compatibility Axiom, Small Cancellation Axiom {\rm(}with respect to $\tau+1$ small pieces{\rm)} and at least one of Isolation Axioms.
\end{defin}

In the further argument we assume that $\tau \geqslant 10$ (recall that in a small cancellation group we require that every relator is a product of not less than $7$ small pieces, see \cite{LS}).

\subsection{Towards a filtration on $\Galg$: multi-turns, replacements, virtual members of the chart and numerical characteristics of monomials }\label{filtration_1}

All the way further we  will study the ring ${\mathcal A} = \Qalg$, with $\Rel$ subject to three small cancellation conditions. 

Let $U$ be a word and $\widehat{U}$ be its subword. We call the triple that consists of $U$, $\widehat{U}$ and the position of $\widehat{U}$ in $U$ \emph{an occurrence of $\widehat{U}$ in $U$}. In fact, we consider occurrences of the form $a \in \Mon$ in $U$, that is $U = LaR$, where $L$, $R$ can be empty. Since $a \in \Mon$, there exists a polynomial $p \in \Rel$ such that $a$ is a monomial of $p$. An {\it overlap} is defined as a common part of two occurrences. Under maximal occurrence we mean an occurrence of a  monomial of $\Mon$ which is not contained in a bigger such occurrence. We shall underline that the a common part of two maximal occurrences is a small piece.



Now we indicate a ring-theoretic counterpart of the notion  of   turn.

\label{chart_def}
\begin{defin}
Let $U$ be a monomial. We define {the chart of $U$} as the set of all maximal occurrences of monomials of $\Mon$ in $U$. The maximal occurrences $m_i \in \Mon$ in $U$ such that $\SPM(m_i) \geqslant \tau$ are called {members of the chart}.
\end{defin}
This means that we count as members of the chart only big  occurrences of monomials from $\Mon$. 
Now we define a multi-turn that  is a ring-theoretic analog of a group turn. 

In the case of groups  we have the following situation.
Let $G$ be a small cancellation group, $R_i = M_1M_2^{-1}$ be a relator of its small cancellation presentation. Assume $LM_1R$ and $LM_2R$ are two words, then the transition from $LM_1R$ to $LM_2R$
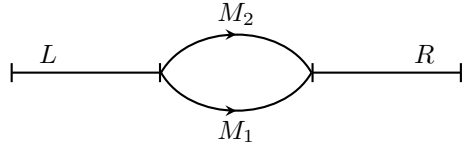
\vspace{0.1cm}

\begin{center}
\begin{tikzpicture}
\draw[|-|, black, thick] (0,0)--(2,0) node [near start, above] {$L$};
\draw[black, thick, arrow=0.5] (2,0) to [bend left=60] node [above] {$M_2$} (4,0);
\draw[black, thick, arrow=0.5] (2,0) to [bend right=60] node [below] {$M_1$} (4,0);
\draw[|-|, black, thick] (4,0)--(6,0) node [near end, above] {$R$};
\end{tikzpicture}
\vspace{0.1cm}
\end{center}
is called {\it a turn} of an occurrence of the subrelation $M_1$ (to its complement $M_2$).
Analogously, in our case we define a multi-turn.



\begin{defin}
\label{multiturn_def}
Let $p = \sum\limits_{j = 1}^{n} \alpha_{j} a_{j} \in \Rel$. For every $h = 1, \ldots, n$ we call the transition
$$
a_h\longmapsto \sum\limits_{{j = 1, \ \\ j\neq h}}^{n} (-\alpha_{h}^{-1} \alpha_{j}a_{j}),
$$
an {elementary multi-turn of $a_h$ with respect to $p$}.

Let $p = \sum\limits_{j = 1}^{n} \alpha_{j} a_{j} \in \Rel$. Let $a_h$ be a maximal occurrence in $U$, $U = La_hR$. The transformation
$$
\label{multi_turn_def}
U = La_hR \longmapsto \sum\limits_{{j = 0, \  \\ j\neq h}}^{n} (-\alpha_{h}^{-1} \alpha_{j}La_{j}R)
$$
with the further cancellations if there are any, is called {a multi-turn of the occurrence $a_h$ in $U$ that comes from an elementary multi-turn $a_h \mapsto \sum_{{j = 1, \ \\ j\neq h}}^{n}(-\alpha_h^{-1}\alpha_j a_j)$}. Obviously,
$$
U - \sum\limits_{{j = 0, \ \\ j\neq h}}^{n} (-\alpha_{h}^{-1} \alpha_{j}La_{j}R) =  \alpha_{h}^{-1}LpR \in \Ideal.
$$
In this case the polynomial $LpR = \sum_{j = 1}^n\alpha_j La_jR$ {\rm (}after the cancellations{\rm )} is called {\it a layout of the multi-turn}.
\end{defin}

In what follows we undertake a very detailed study of the influence
of multi-turns on  charts of the
monomials.  We will  trace transformation of a chart under the given multi-turn or set of multi-turns. 
We also take care of transformations of individual monomials $U_h=La_hR \mapsto U_j=La_jR$ called {\it replacements}.

 Applying the multi-turns of $a_h$ in $U_h = La_hR$ we arrive at monomials  $U_j = La_jR$. We describe precisely {\it how the corresponding
maximal occurrences in $U_j$ look like comparatively to maximal occurrences in $U_h$                                    }.

We consider three variants for the resulting monomial $U_j=La_jR$: 
 $a_j$ is not a small piece; $a_j$ is  a small piece; $a_j$ is 1.
  We show that in the first case the structure of the chart remains almost stable after a multi-turn,
   in the second case the replacement $a_h$ by $a_j$ can cause merging and  restructuring of the chart, and in the third case strong cancellations resulting in complete modification of the chart are possible.

We produce the full list of all appearing arrangements of maximal occurrences.
The calculations are based on thorough analysis of all combinatorial possibilities. This list is in fact a Theorem that provides ground to further considerations towards a filtration on $\Galg$.

{\it Our goal is
 constructing
a special ordering on monomials. This ordering is far from being usual $DegLex$-order. In more precise terms our objective is to build numerical characteristic of a chart that allows to define a filtration on monomials which behaves well with respect to replacements of the monomials caused by multi-turns.}

On the way we have to treat several caveats. When we define members of
a chart in the terms of their $\Lambda$-measure, such definition is not stable enough
under multi-turns. So, we define a quite delicate notion of a {\it virtual
member} of a chart. Virtual members of the chart are those occurrences $b$ which  originally are not necessarily  members of the chart but they are rather big with $\Lambda(b)\geqslant \tau-2$,  and after a series of {\it admissible transformations} become members of the chart. In turn, {\it  admissible replacements} are those  $a_h\mapsto a_j$  that  take sufficiently long monomials $a_h\in U_h$ with $\Lambda(a_h)\geqslant \tau - 2$ to monomials $a_j\in U_j$ which are not fully covered by images of elements of   $\longmo{U_h}  \setminus \lbrace a_h\rbrace$ in $U_j$. Here $\longmo{U_h}$ stands for the set of all  maximal occurrences in $U_h$ of $\SPM$-measure $\geqslant 3$.


Let $U$ be a monomial. Consider subsets of $\mo{U}$ that cover the same letters in $U$ as the whole $\mo{U}$. A covering of such type consisting of the smallest number of elements is called {\it a minimal covering}.  Of course, such covering is not, necessarily, unique.

Given a monomial $U$, we define $\mincov{U}$
to be the number of elements in a minimal covering of $U$.  Denote the number of virtual members of the chart of $U$ by  $\nvirt{U}$. It is clear that $\nvirt{U}\leqslant \mincov{U}$.


The next proposition aggregates all calculations beforehand. 
\begin{proposition}
\label{non_increasing_parameter}
Assume $U_h$ is a monomial, $a_h$ is a virtual member of the chart of $U_h$. Let $a_h$ and $a_j$ be incident monomials. Consider the replacement $a_h \mapsto a_j$ in $U_h$. Let $U_j$ be the resulting monomial. If $a_j$ is a virtual member of the chart of $U_j$, then $\mincov{U_h} = \mincov{U_j}$ and $\nvirt{U_h} = \nvirt{U_j}$. If $a_j$ is not a virtual member of the chart of $U_j$, then either $\mincov{U_j} < \mincov{U_h}$, or $\mincov{U_j} = \mincov{U_h}$ but $\nvirt{U_j} < \nvirt{U_h}$.
\end{proposition}

\begin{defin}
Let $U$ be a monomial. We introduce $f$-characteristic  $U$ by the rule:
$$
\label{f_char}
f(U) = (\mincov{U}, \nvirt{U)}).
$$
If $U_1$ and $U_2$ are monomials, we say that $f(U_1) < f(U_2)$ if and only if either $\mincov{U_1} < \mincov{U_2}$, or $\mincov{U_1} = \mincov{U_2}$ but $\nvirt{U_1} < \nvirt{U_2}$.
\end{defin}

We define {\it derived monomials} of $U$ as the result of applying of a sequence of replacements of virtual members of the chart by incident monomials, starting from $U$.

%
\begin{lemma}
\label{estimation_value_property}
Assume $U$ and $Z$ are monomials, $Z$ is a derived monomial of $U$. Then $f(Z) \leqslant f(U)$. Moreover, $f(Z) < f(U)$ if and only if in the corresponding sequence of replacements there exists at least one replacement of the form $La_hR \mapsto La_jR$ such that $a_h$ is a virtual member of the chart of $La_hR$ and $a_j$ is not a virtual member of the chart of $La_jR$.
\end{lemma}

The introduced    $f$-characteristic gives rise to a certain function $t$ on natural numbers defined as follows. We put $t(0) = (0, 0)$. Assume $t(n) = (r, s)$, then we put

$$
t(n + 1)=\left\lbrace \begin{array}{ccc} (r,&s+1)&\textit{ if } r > s, \\ (r+1,&0)&\textit{ if } r = s.  \\ \end{array}\right.
$$

\begin{defin} We define an increasing filtration on $\Galg$ by the rule:
\label{filtration_def}
$$
\Ft_n(\Galg) = \langle Z \mid Z \in \Fr, f(Z) \leqslant t(n)\rangle.
$$
\end{defin}
That is, the space $\Ft_n(\Galg)$ is generated by all monomials with $f$-characteristics not greater than $t(n)$.


\subsection{Derived monomials and dependencies}\label{filtration_2}

We need a set of new notions. Let $U$ be a monomial. By $\langle U\rangle_d$ we denote a linear subspace of $\Ft_n(\Galg)$ generated by all derived monomials of $U$. By $\Low\DSpace{U}$ we denote the subspace generated by all derived monomials of $U$ with $f$-characteristic smaller than $f(U)$. The next principal object is the set of {\it dependencies}, defined as follogradingws. Suppose $Y$ is a subspace of $\Galg$ linearly generated by a set of monomials and closed under taking derived monomials. We take the set of all layouts of multi-turns of virtual members of the chart of monomials of  $Y$ and look at its linear envelope $\Dp(Y)$, which is our set of dependencies related to $Y$. We prove that $\Dp(\Galg) = \Ideal$.

The  key statement is the following Proposition which describes  nice interaction between dependencies and filtration:
\begin{proposition}\label{mainprop}
$$
\Dp(\Ft_n(\Galg)) \cap \Ft_{n - 1}(\Galg) = \Dp(\Ft_{n - 1}(\Galg)).
$$
\end{proposition}
\noindent
This proposition yields
\begin{proposition}
\label{fall_to_smaller_subspace}
Suppose $X, Y$ are subspaces of $\Galg$ generated by monomials and closed under taking derived monomials, $Y \subseteq X$. Then $\Dp(X)\cap Y = \Dp(Y)$.
\end{proposition}
Proof of Proposition \ref{mainprop} is based on Main Lemma. Namely,

\begin{lemma}[\textbf{Main Lemma}]\label{mainlemma}
Let $U$ be an arbitrary monomial, $U \in \Ft_n(\Galg) \setminus \Ft_{n - 1}(\Galg)$. Then
$$
\Dp\DSpace{U} \cap \Low\DSpace{U} \subseteq \Dp(\Ft_{n - 1}(\Galg)).
$$
\end{lemma}

Here is the place to make  some comments. Main Lemma says that there is a natural interaction between dependencies and   reduction of $f$-characteristic,  and this interaction causes descending in the filtration. 
This yields, in essence, 
that in the quotient algebra $\Qalg$ {\it there are no unexpected linear dependencies}. But, first, one has to explain what are the {\it expected} linear dependencies.

Consider the filtration $\Ft_n (\Galg)$, $n \geqslant 0$, on $\Galg$ defined as above.
Let $U \in \Ft_n (\Galg)$ be a monomial such that its chart has $m$ virtual members $u^{(i)}$, $U = L^{(i)} u^{(i)} R^{(i)}$, $ i= 1, 2, \ldots , m$. For any $p \in \Rel$ of the form $p = \alpha u^{(i)} + \sum_{j=1}^k \alpha_j a_j$, $\alpha \neq 0$, we consider the polynomial $L^{(i)} p  R^{(i)} \in \Galg$. All such polynomials obviously belong to $\Ft_n (\Galg) \cap \Ideal$ and regarded as expected dependencies.  We shall emphasize that in case the relations $\Rel$ do not satisfy special conditions, the term $\Ft_n (\Galg) \cap \Ideal$ may contain also arbitrary unexpected dependencies.

In fact,  Proposition \ref{fall_to_smaller_subspace} claims that the opposite is also true. In more detail, Proposition \ref{fall_to_smaller_subspace} implies that $\Ft_n (\Galg) \cap \Ideal = \Ft_n(\Galg) \cap \Dp(\Galg) = \Dp(\Ft_n(\Galg))$. That is, $\Ft_n(\Galg) \cap \Ideal$ is linearly generated by expected linear dependencies related to $\Ft_n(\Galg)$. This can be restated as follows.


\begin{theorem}\label{mainlemmath}
$\Ft_n (\Galg) \cap \Ideal$
is linearly spanned by all  polynomials of the form $L^{(i)} p  R^{(i)}$, $i=1,  \ldots , m,$ for all monomials $U \in \Ft_n (\Galg)$ and polynomials $p \in \Rel$ as above, $n \geqslant 0$.
\end{theorem}

\subsection{Grading on small cancellation ring}\label{grading}

%
First of all,  it can be seen  that  $\Dp(\Galg) = \Ideal$. The quotient space $\Qalg$ naturally inherits the filtration from $\Galg$, namely,
$$
\Ft_n (\Qalg) = (\Ft_n(\Galg) + \Dp(\Galg)) / \Dp(\Galg) = (\Ft_n(\Galg) + \Ideal) / \Ideal.
$$

We define a grading on $\Qalg$ by the rule:
$$
\Gr(\Qalg) = \bigoplus\limits_{n = 0}^{\infty}\Gr_n(\Qalg)=\bigoplus\limits_{n = 0}^{\infty}\Ft_n (\Qalg) / \Ft_{n - 1} (\Qalg).
$$
The next theorem establishes the compatibility of the filtration and the corresponding grading on~$\Galg$ with the space of dependencies~$\Dp(\Galg)$. It states that
\begin{theorem}
\label{structure_of_quotient_space}
$$\Gr_n(\Qalg) \cong \Ft_n(\Galg) / (\Dp(\Ft_n(\Galg)) + \Ft_{n - 1}(\Galg)).$$
\end{theorem}

\subsection{Non-triviality of $\Qalg$. Construction of a basis of $\Qalg$}\label{basis_nontriv}
\begin{lemma}
\label{non_trivial_spaces_existence}
Let $\lbrace V_i\rbrace_{i \in I}$ be all different spaces $\lbrace\DSpace{Z} \mid Z \in \Fr\rbrace$. Then not all spaces $V_i / (\Dp(V_i) + \Low(V_i))$, $i\in I$, are trivial. Namely, the space $\DSpace{X} / (\Dp\DSpace{X} + \Low\DSpace{X})$, where $X$ is a monomial with no virtual members of the chart, is always non-trivial, and of dimension~$1$. In particular, $\DSpace{1} / (\Dp\DSpace{1} + \Low\DSpace{1}) \neq 0$, where $1$ is the empty word.
\end{lemma}

Proof. Let $X$ be a monomial with no virtual members of the chart. Then there are no derived monomials of $X$ except $X$ itself, and there are no multi-turns of virtual members of the chart of $X$. So, by definition, $\DSpace{X}$ is linearly generated by $X$ and, therefore, is of dimension~$1$; $\Dp\DSpace{X} = 0$; $\Low\DSpace{X} = 0$. Therefore,
$$
\DSpace{X} / (\Dp\DSpace{X} + \Low\DSpace{X}) = \DSpace{X} = \langle X\rangle \neq 0,
$$
and $\DSpace{X} / (\Dp\DSpace{X} + \Low\DSpace{X}) $ is of dimension~$1$.

By definition, the empty word~$1$ is a small piece. Therefore, $1$ has no virtual members of the chart. So, it follows from the above that $\DSpace{1} / (\Dp\DSpace{1} + \Low\DSpace{1}) \neq 0$.

Now we can prove that the quotient ring $\Qalg$ is non-trivial.
\begin{corollary}
\label{non_trivial_quotient}
The quotient ring $\Qalg$ is non-trivial.
\end{corollary}
Proof.
Let $U$ be a monomial. Consider the space $\DSpace{U}$ and the corresponding subspace in $\Qalg$, namely, $(\DSpace{U} + \Ideal) / \Ideal$. From the isomorphism theorem it follows that
$$
(\DSpace{U} + \Ideal) / \Ideal \cong \DSpace{U} / (\DSpace{U}\cap \Ideal).
$$
Recall that $\Ideal = \Dp(\Galg)$. From Proposition~\ref{fall_to_smaller_subspace} it follows that $\DSpace{U} \cap \Dp(\Galg) = \Dp\DSpace{U}$. Hence,
$$
(\DSpace{U} + \Ideal) / \Ideal \cong \DSpace{U} / \Dp\DSpace{U}.
$$

By Lemma~\ref{non_trivial_spaces_existence}, there exists a space $\DSpace{U_0}$, $U_0 \in \Fr$, such that $\DSpace{U_0} / (\Dp\DSpace{U_0} + \Low\DSpace{U_0}) \neq 0$. Hence, we see that $\DSpace{U_0} / \Dp\DSpace{U_0} \neq 0$ and $(\DSpace{U_0} + \Ideal) / \Ideal \neq 0$. So, there exists a non-trivial subspace of $\Qalg$. Thus, $\Qalg$ itself is non-trivial.

Now we are able, at last, to describe a basis of $\Qalg$. This is done in two steps. First, we construct  a basis for non-trivial graded components of our filtration on $\Qalg$:
$$
\Gr_n(\Qalg) = \Ft_n(\Qalg) / \Ft_{n - 1}(\Qalg).
$$
Given $n$ we consider the set of spaces $\lbrace\DSpace{Z} \mid Z \in \Fr, Z \in \Ft_n(\Galg)  \setminus \Ft_{n - 1}(\Galg)\rbrace $, such that $\DSpace{Z} / (\Dp\DSpace{Z} + \Low\DSpace{Z}) \neq 0$. Let $\lbrace V_{i}^{(n)}\rbrace_{i \in I^{(n)}}$ be all different spaces from this set. Then,
$$
\Gr_n(\Qalg) \cong \bigoplus\limits_{i \in I^{(n)}} V_{i}^{(n)} /(\Dp(V_{i}^{(n)}) + \Low(V_{i}^{(n)})).
$$
Assume $\lbrace \overline{W}^{(i, n)}_j\rbrace_j$ is a basis of $V_{i}^{(n)}/(\Dp(V_{i}^{(n)}) + \Low(V_{i}^{(n)}))$, $i \in I^{(n)}$. Let $W^{(i, n)}_j \in V_{i}^{(n)}$ be an arbitrary representative of the coset $\overline{W}^{(i, n)}_j$. Then
$$
\bigcup\limits_{i \in I^{(n)}}\left\lbrace W^{(i, n)}_j + \Ideal + \Ft_{n - 1}(\Qalg)\right\rbrace_j
$$
is a basis of $\Gr_n(\Qalg)$.


Finally, the next Theorem describes a basis of $\Qalg$. We have
\begin{theorem}\label{mainth} Let $\lbrace V_i\rbrace_{i \in I}$ be all different spaces $\lbrace\DSpace{Z} \mid Z \in \Fr\rbrace$. Then
$$
\Qalg \cong \bigoplus\limits_{i \in I} V_i/(\Dp(V_i) + \Low(V_i)),
$$
as vector spaces, and the right-hand side is explicitly described via a tensor product of subspaces. 

Assume $\lbrace \overline{W}^{(i)}_j \rbrace_j$ is a basis of $V_i/(\Dp(V_i) + \Low(V_i))$, $i \in I$. Let $W^{(i)}_j \in V_i$ be an arbitrary representative of the coset $\overline{W}^{(i)}_j$. Then
$$
\bigcup\limits_{i \in I}\left\lbrace W^{(i)}_j  + \Ideal\right\rbrace_j
$$
is a basis of $\Qalg$.
\end{theorem}

\subsection{Examples, algorithmic properties 
}\label{examples_algorithms}

We study algorithmic properties of the constructed small cancellation  ring. We show that they are as expected to be for small cancellation objects and similar in a sense to the ones valid for small cancellation groups. However, in the ring case  the essential peculiarities arise in many places. Recall that small cancellation groups enjoy Dehn's  algorithm \cite{LS}. 
In this section we define and study  a corresponding greedy  algorithm for rings.

Let  a ring ${\mathcal A} = \Qalg$ with small cancellation condition be given. We extend  our set of relations $\Rel$ to a certain  additive closure $\Add(\Rel)$. It is important that $\Rel = \Add(\Rel)$ for natural examples of the ring ${\mathcal A}$ considered below. We define a linear order on all monomials, based on $f$-characteristic and other  considerations,  and denote it by $<_f$. Then, given the  order $<_f$ and the set $\Add(\Rel)$, we define a special greedy algorithm  (with external source of knowledge) for small cancellation rings. This algorithm has the similar meaning as {\it Dehn's algorithm} does for the case of groups.  Denote it by $\GreedyAlg(<_f, \Add(\Rel))$.

Recall that given a small cancellation group $G = \langle \mathcal{X} \mid \Rel_G\rangle$, a word $W$ from a free group
is equal to $1$ in $G$ if and only if Dehn's algorithm, starting from $W$, terminates at $1$, \cite{LS}. Our Theorem \ref{alg}
establishes the similar properties in much more complicated situation of rings.

Namely, assume $W_1, \ldots, W_k$ are different monomials. We take an element $\sum_{i = 1}^k \gamma_i W_i \in \Galg$, $\gamma_i \neq 0$.
\begin{theorem}\label{alg}
The following statements are equivalent:
\begin{itemize}
\item [$\bullet$]
\label{possibility_of_reduction1}
some branch of the algorithm $\GreedyAlg(<_f, \Add(\Rel))$, starting from $\sum_{i = 1}^k \gamma_i W_i$, terminates at~$0$;
\item [$\bullet$]
\label{possibility_of_reduction2}
$\sum_{i = 1}^k \gamma_i W_i \in \Ideal$;
\item [$\bullet$]
\label{possibility_of_reduction3}
every branch of the algorithm $\GreedyAlg(<_f, \Add(\Rel))$, starting from $\sum_{i = 1}^k \gamma_i W_i$, terminates at~$0$.
\end{itemize}
\end{theorem}
\begin{corollary}\label{mainth_algorithms} We have
\begin{itemize}
\item [$\bullet$]
 $\GreedyAlg(<_f, \Add(\Rel))$ solves the Ideal Membership Problem for~$\Ideal$,
\item [$\bullet$]
$\Add(\Rel)$ is a Gr\"{o}bner basis of the ideal $\Ideal$ with respect to monomial ordering $<_f$.
\end{itemize}
\end{corollary}

We give two examples of small cancellation rings. One can check that the group algebra of a small cancellation group satisfying   a small cancellation condition with $C(m)$ for $m\geqslant 22$ (see \cite{LS}) is a small cancellation ring. Another example is a ring constructed in \cite{AKPR}. This is  a quotient ring $\mathbb{Z}_2\Fr / \Ideal$, where $\mathbb{Z}_2\Fr$ is the group algebra of a free group $\Fr$ over the field $\mathbb{Z}_2$, and the ideal $\Ideal$ is generated by a single trinomial $1 + v + vw$, where $v$ is a complicated word depending on $w$.
The ring $\mathbb{Z}_2\Fr / \Ideal$ is of special interest, since $(1 + w)^{-1} = v$ in it. Thus, binomial  $1 + w$ becomes invertible.


{\bf Acknowledgements.} The research of the first, second and the third authors was supported by  ISF grant 1994/20 and the Emmy Noether Research Institute for Mathematics. The research of the first author was also supported by ISF fellowship. The research of the second author was also supported by the Russian science foundation, grant  17-11-01377.

We are very grateful to I.Kapovich, B.Kunyavskii and D.Osin   for invaluable cooperation.


\begin{thebibliography}{00}

\bibitem{A}[A]
	S.~I.~Adian, {\it The Burnside problem and identities in groups,} Nauka, Moscow, 1975 , 335 pp.

\bibitem{AKPR}[AKPR]
A.~Atkarskaya, A.~Kanel-Belov, E.~Plotkin, E.~Rips, {\it Construction of a quotient ring of $Z_2F$ in which a binomial 1+w is invertible using small cancellation methods},
in "Groups, Algebras, and Identities" AMS, Contemporary Mathematics,
{\bf 726} (2019), Israel Mathematical Conferences Proceedings,  1--76.

\bibitem{AKPR1}[AKPR1]
A.~Atkarskaya, A.~Kanel-Belov, E.~Plotkin, E.~Rips, {\it Group-like Small Cancellation Theory for Rings}, Arxiv, 244pp.










\bibitem{Bo}[Bo]
B.~Bowditch.  {\it A course on geometric group theory},  MSJ Memoirs. 16. Tokyo: Mathematical Society of Japan.





\bibitem{DK}[DK]
C.~Drutu, M.~Kapovich, {\it Geometric group theory}, Colloquiunm Publications, AMS, {\bf 63}, 2018, 807pp.


\bibitem{Gr1}[Gr1]
M.~Gromov, {\it Infinite groups as geometric objects},
Proc. Int. Congress Math., Warsaw,
1983, Amer. Math. Soc.,  {\bf 1} (1984), 385--392.


\bibitem{Gr2}[Gr2]
M.~Gromov, {\it Hyperbolic Groups},
 "Essays in Group Theory" (G. M. Gersten, ed.),
{\bf 8} (1987),  MSRI Publ., Springer, New York, 75--263.

\bibitem{Gu}[Gu]
V.~Guba, {\it Finitely generated complete groups}, Izv. Akad. Nauk SSSR Ser. Mat. {\bf 50} (1986), 883-924.

\bibitem{GS}[GS]
V.~Guba, M.~Sapir, {\it Diagram Groups}, Memoirs of the American Mathematical Society
(1997), 117 pp.







\bibitem{Hig}[Hig]
P.~M.~Higgins, {\it Techniques of semigroup theory}, Oxford University Press,
Oxford, (1992).

\bibitem{Ivanov}[Ivanov]
S.~Ivanov, {\it The free Burnside groups of sufficiently large exponents}, Internat. J. Algebra Comput. {\bf 4} (1994), no. 1-2, ii+308pp.



 \bibitem{Ly}[Ly]
 R.~Lyndon, {\it On Dehn's algorithm}, Math. Ann. {\bf 166}, (1966), 208--228.


 \bibitem{LS}[LS]
R.~Lyndon,  P.~Schupp, {\it Combinatorial group theory. Reprint of the 1977 edition},
Classics in Mathematics. Springer-Verlag, Berlin.
(2001).

\bibitem{Lys}[Lys]
I.~Lysenok (1996). {\it Infinite Burnside groups of even exponent}. Izv. Math. {\bf 60:3} (1996),  453--654.



\bibitem{NA1}[NA1]
P.S.~Novikov, S.I.~Adian, {\it Infinite periodic groups. I},
Izvestia Akademii Nauk SSSR. Ser. Mat.,
{\bf 32} (1968), no. 1, 212--244.

\bibitem{NA2}[NA2]
P.S.~Novikov, S.I.~Adian, {\it Infinite periodic groups. II},
Izvestia Akademii Nauk SSSR. Ser. Mat.,
{\bf 32} (1968), no. 2, 251--524.

\bibitem{NA3}[NA3]
P.S.~Novikov, S.I.~Adian, {\it Infinite periodic groups. III},
Izvestia Akademii Nauk SSSR. Ser. Mat.,
{\bf 32} (1968), no. 3, 709--731.

\bibitem{Ol1}[Ol1]
A.~Olshanskii, {\it Geometry of defining relations in groups. Translated from the 1989 Russian original by Yu. A. Bakhturin},
Mathematics and its Applications (Soviet Series), Kluwer Academic Publishers Group, Dordrecht,
{\bf 70} (1991).

\bibitem{Ol2}[Ol2]
A.~Olshanskii, {\it An infinite  group with subgroups of prime orders},
Math. USSR Izv.
{\bf 16} (1981), 279-289;
translation of Izvestia Akad. Nauk SSSR Ser. Matem. {\bf 44} (1980), 309--321.

\bibitem{Ol3}[Ol3]
A.~Olshanskii, {\it Groups of bounded period with subgroups of prime order},
Algebra and Logic,
{\bf 21} (1983), 369--418; Translation of Algebra i Logika {\bf 21} (1982), 553--618.

\bibitem{Rips}[Rips]
E.Rips,  {\it Generalized small cancellation theory and applications I},  Israel J. Math., {\bf 41} (1982), 1--146.



\bibitem{S}[S]
M.~Sapir, Combinatorial algebra: syntax, Springer Monographs in Mathematics, Springer, Cham, 2014, 355 pp.






\end{thebibliography}
\end{document}